\magnification1200
\input epsf

\centerline{\bf Hyperbolic Coxeter groups with  Sierpi\'nski
carpet  boundary}

\medskip
\centerline{{Jacek \'Swi\c atkowski}
\footnote{*}{This work was partially supported by the Polish National Science Centre (NCN), grant 2012/06/A/ST1/00259.}}

\medskip
\centerline{Instytut Matematyczny, Uniwersytet Wroc\l awski}
\centerline{pl. Grunwaldzki 2/4, 50-384 Wroc\l aw, Poland}
\centerline{e-mail: {\tt swiatkow@math.uni.wroc.pl}}
%\medskip
%\centerline{February 4, 2015}

\bigskip\noindent
{\bf Abstract.}
We give a necessary and sufficient condition for a hyperbolic Coxeter group
with planar nerve to have Sierpi\'nski curve as its Gromov boundary.

\bigskip\noindent
{\bf Mathematics Subjest Classification (2010).} 20F67; 20F55, 20F65.


\bigskip

\noindent
{\bf Introduction.} 

\medskip
In the survey [KB] concerning boundaries of hyperbolic groups
the authors wrote: 
{\it classifying groups which have the Sierpi\'nski carpet as the boundary
remains a challenging open problem}. In this paper we deal with this problem 
for hyperbolic Coxeter groups.
We first state (as Theorem 1 below)
our main result in the framework of right angled hyperbolic
Coxeter groups, as this does not require much preparations,
and is particularly simple. Then, after introducing appropriate
terminology, we present the main result in its full generality
(as Theorem 2 below).

\medskip
%\noindent
%{\bf Definition.}
We call
a simplicial complex {\it unseparable} if it is connected, has no separating
simplex, no separating pair of nonadjacent vertices, and no separating
full subcomplex isomorphic to the (simplicial) suspension of a simplex.

\medskip\noindent
{\bf Theorem 1.}
{\it Let $(W,S)$ be a right angled Coxeter system such that the group $W$ is word hyperbolic and the nerve
$L$ of the system is planar, distinct from a simplex and from a triangulation of the 2-sphere $S^2$. 
Then the Gromov boundary $\partial W$
is homeomorphic to the Sierpi\'nski curve if and only if $L$
is unseparable.}

\centerline{\epsfbox{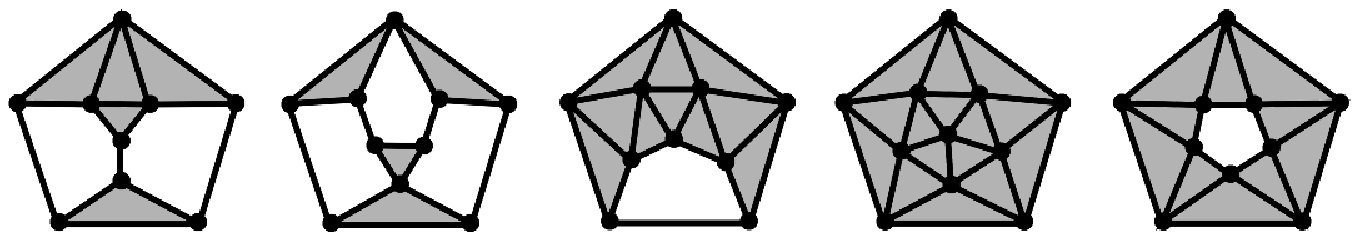}}

\medskip
\centerline{Figure 1. Examples of nerves of hyperbolic right-angled Coxeter groups}
\centerline{with Sierpi\'nski carpet boundaries.}

\bigskip

%\medskip
Recall that for a right-angled Coxeter system $(W,S)$ the group $W$ is word hyperbolic
iff the nerve $L$ of this system (which is a flag simplicial complex) satisfies {\it no empty square} condition, i.e. any polygonal cycle of length 4 in the 1-skeleton of $L$ has at least one diagonal (see the comment at the bottom of p. 233 in [Da]).
Examples of some nerves of hyperbolic right-angled Coxeter groups with Sierpi\'nski carpet boundaries are presented at Figure 1.

\smallskip
We extend the concept of unseparabilty to appropriately understood nerves of arbitrary (not only right-angled) Coxeter systems.

%\break
\medskip\noindent
{\bf Definition.}
\item{(1)}
A {\it labelled nerve} $L^\bullet$ 
of a Coxeter system $(W,S)$ is the nerve
$L$ of
$(W,S)$ equipped with the labelling, which to any edge of $L$
associates the corresponding entry from the Coxeter matrix
of this system. More precisely, if $e$ is an edge of $L$ with
vertices $s,t\in S$ then the label associated to $e$ is equal to
the exponent $m_{st}$ in the relator $(st)^{m_{st}}$ 
from the standard presentation
for $W$ corresponding to the system $(W,S)$.

\item{(2)}
A {\it labelled suspension} in a labelled nerve $L^\bullet$ is a full subcomplex
$\Lambda$ of $L$ isomorphic to the simplicial suspension of a simplex,
$\Lambda=\{ t,s \}*\sigma$, such that any edge in $\Lambda$
adjacent to $t$ or $s$ has label 2.

\item{(3)}
The labelled nerve $L^\bullet$ of a Coxeter system is {\it unseparable}
if it is connected, has no separating
simplex, no separating pair of nonadjacent vertices, and no separating
labelled suspension.

\medskip\noindent
{\bf Remark.}
Unseparability of the labelled nerve of $(W,S)$ is equivalent to the fact
that $W$ has no visual splitting 
(in the sense of Mihalik and Tschantz [MT])
along a finite or a 2-ended subgroup.

\medskip
Theorem 1 is a special case of the following more general result.
In its statement, the condition that $L^\bullet$ is a {\it labelled wheel}
means that $L$ is isomorphic to the simplicial cone over a triangulation
of $S^1$, and the labels at all edges adjacent to the cone vertex
are equal to 2. (If $L^\bullet$ is a labelled wheel then $W$ is easily
seen to be a virtually surface group, so its boundary is then
a circle.)

\medskip\noindent
{\bf Theorem 2.}
{\it Let $(W,S)$ be a Coxeter system such that the group $W$ is word hyperbolic and
the nerve of the system is planar,
distinct from a simplex and from a triangulation of $S^1$ or $S^2$. 
Then the Gromov boundary $\partial W$
is homeomorphic to the Sierpinski curve if and only if 
the labelled nerve $L^\bullet$ of $(W,S)$
is distinct from a labelled wheel and unseparable.}

\medskip
Recall that a Coxeter group is word hyperbolic iff it contains no affine special subgroup
of rank $\ge3$, and no special subgroup being the direct sum of two infinite special 
subgroups (see [Da], Corollary 12.6.3 on p. 241).
Figure 2 presents a sample of labelled nerves of not right-angled hyperbolic
Coxeter groups with Sierpi\'nski carpet boundaries. In this figure we use the convention
that an edge with no label indicated is labelled with 2.

\medskip\noindent
{\bf Remarks.}
\item{(1)}
The following two conditions from the statement of Theorem 2
do not have their explicitely
expressed counterparts in the statement of Theorem 1: 
\itemitem{$\bullet$} that $L$ is not a triangulation of $S^1$, and
\itemitem{$\bullet$} that $L^\bullet$ is not a labelled wheel.

\item{}
However, the counterparts of these conditions follow
automatically from the remaining assumptions of Theorem 1.
More precisely, $L$ cannot be the triangulation of $S^1$
consisting of 3 edges (the 3-cycle), or the wheel equal to the simplicial cone over the 3-cycle (the 3-wheel), because nerves of right angled Coxeter systems are flag. Furthermore, $L$ cannot be any other triangulation of $S^1$
because of the assumption that it has no separating pair of nonadjacent
vertices. Similarly, $L$ cannot be any other wheel because of the
assumption that it has no separating suspension.
This justifies that Theorem 2 is an extension of Theorem 1.
%\item{(2)} As a consequence of the above remark (1),
%in the statement of Theorem 2 one can replace
%the assumption that $L$ is not a triangulation of $S^1$
%with the weaker assumption that it is not the 3-cycle.
%Similarly, the condition that $L^\bullet$ is not a labelled wheel
%may be replaced with the weaker condition that it is not
%a labelled 3-wheel. 
%[TO CHYBA PRZEGADANE I MOZNA POMINAC]

\item{(2)}
%\medskip
Note that Theorem 2, among others, determines all hyperbolic Coxeter groups
with Gromov boundary homeomorphic to the Sierpi\'nski curve, 
and whose nerves are
\itemitem{$\bullet$} planar graphs,
\itemitem{$\bullet$} triangulations of the 2-disk,
\itemitem{$\bullet$} triangulations of other planar surfaces.

\item{}
%\noindent 
It seems that triangulations of surfaces were never previously considered
as potential candidates for nerves of Coxeter groups with the Sierpi\'nski
carpet boundary.
Some limited results concerning planar graphs as nerves of
hyperbolic Coxeter groups with carpet boundaries were obtained
by N. Benakli, as well as by N. Chinen and T. Hosaka (these works 
seem not
to be yet fully documented).
%This significantly extends all previously known results concerning
%Coxeter groups with Sierpi\'nski curve boundaries.

\item{(3)}
It is possible that, up to product with a finite Coxeter group,
planarity of the nerve is a necessary condition for any Coxeter group
to have the Sierpi\'nski carpet boundary (and even more generally,
to have planar boundary). If this were the case, Theorem 2 would
provide the complete necessary and sufficient condition
for a hyperbolic Coxeter system to have the Sierpi\'nski carpet boundary.
However, we do not resolve this issue in the present paper.

\centerline{\epsfbox{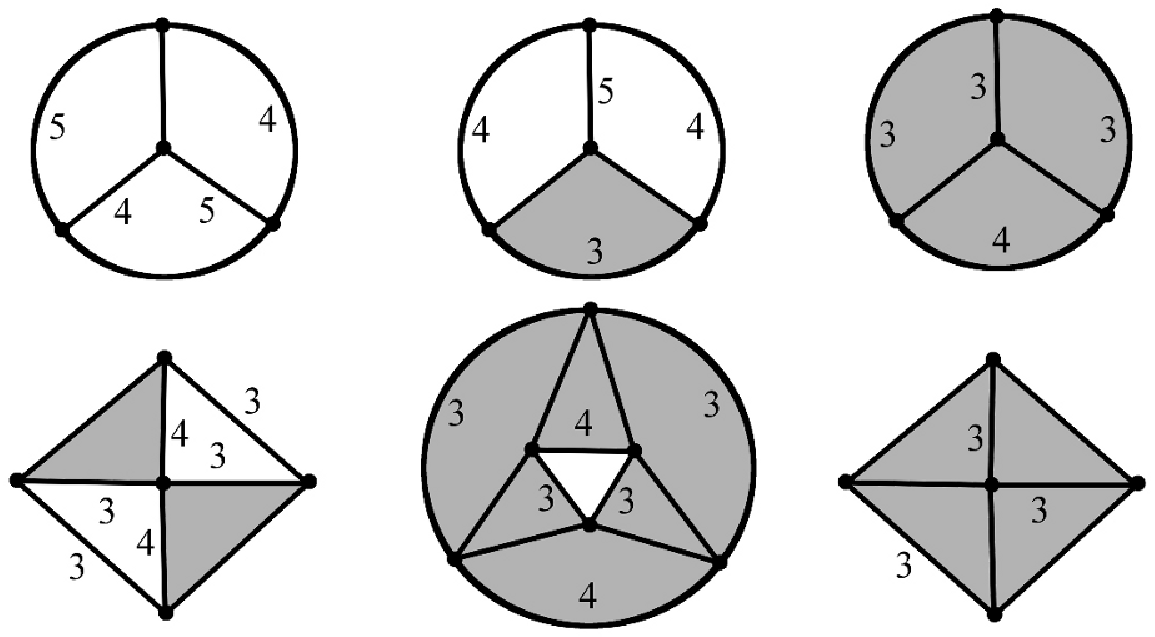}}

\medskip
\centerline{Figure 2. Examples of labelled nerves of hyperbolic not right-angled Coxeter groups}
\centerline{with Sierpi\'nski carpet boundaries.}

%\break
\medskip
The paper is organized as follow. In Section 1 we deal with arbitrary (i.e.
not necessarily hyperbolic) Coxeter groups. We state a conjecture concerning
the occurence of Sierpi\'nski carpet as boundary for such groups, which generalizes
Theorem 2. We then show various partial results related to this conjecture
(which are also the steps in our proof of Theorem 2).

In Section 2 we deal with hyperbolic Coxeter groups, providing those arguments
for showing Theorem 2 which rely upon the hyperbolicity assumption.
We also include comments indicating why these arguments cannot be easily
extended to the general (i.e. not necessarily hyperbolic) case.

The proofs in Sections 1 and 2 are relatively short, and they mainly consist
of appropriate references to various so far known results.

In the appendix we present a short proof of the fact 
(which we use in our arguments)
that boundary
of a special subgroup of a Coxeter group is canonically a subspace
in the boundary of this group. It seems that this folklore result,
known to the experts for at least 20 years, does not yet have
a proof explicitely
presented in the literature of the subject.

\bigskip

\noindent
{\bf 1. A conjecture and some partial results for arbitrary Coxeter groups.}

\medskip
%In this section we present a short proof of the restriction of Theorem 2
%to word hyperbolic Coxeter groups. We present also some partial results 
%concerning the general case, and comments indicating why the short proof 
%does not extend fully.

Recall that given an arbitrary Coxeter system $(W,S)$, there is canonically
associated to it a $CAT(0)$ piecewise euclidean complex $\Sigma=\Sigma(W,S)$,
called {\it the Coxeter-Davis complex} of the system. The boundary $\partial(W,S)$
is defined as the visual boundary of this associated Coxeter-Davis complex.
It is known that when the group $W$ is word hyperbolic, the boundary
$\partial(W,S)$ coincides (up to homeomorphism) with the Gromov boundary
$\partial W$. This allows for the following generalization of Theorem 2,
which we formulate as a conjecture.

\medskip\noindent
{\bf 1.0 Conjecture.}
{\it  Let $(W,S)$ be a Coxeter system with planar nerve 
distinct from a simplex and from a triangulation of $S^1$ or $S^2$. 
Then the boundary $\partial(W,S)$
is homeomorphic to the Sierpinski curve if and only if 
the labelled nerve $L^\bullet$ of $(W,S)$
is distinct from a labelled wheel and unseparable.}

\medskip
Recall that  
the Sierpi\'nski curve has the following characterization due to Whyburn [Wh]:
it is the unique up to homeomorphism compact connected and locally connected metric space of topological dimension 1, which is planar (i.e. can be embedded in the 2-sphere $S^2$), 
and which has no local cut point. For short, we will denote the Sierpi\'nski
curve by $\Pi$.

\medskip
Throughout this section $(W,S)$ is a Coxeter system, $L$ denotes its nerve,
and $L^\bullet$ denotes its labelled nerve.
We start with a lemma yielding  
one of the implications in the above conjecture (and thus also in Theorem 2).

\medskip\noindent
{\bf 1.1 Lemma.}
{\it Unseparability of the labelled nerve $L^\bullet$ is a necessary condition
for $\partial(W,S)$ to be homeomorphic to the Sierpi\'nski curve $\Pi$.
More precisely, it is a necessary condition for $\partial(W,S)$ to be connected
and to have no local cut points.}

\medskip\noindent
{\bf Proof:} 
By Theorem 8.7.2 in [Da], if $L$ is not connected, or has a separating
simplex, then $W$ is either 2-ended or has infinitely many ends.
Consequently, $\partial(W,S)$ is not connected, and hence not
homeomorphic to $\Pi$.

Consider the remaining case when $W$ is 1-ended.
If $L$ has a separating pair of nonadjacent vertices,
or $L^\bullet$ has a separating labelled suspension, then $W$
splits over a 2-ended special subgroup (induced by the corresponding
separating subcomplex). Consequently, the boundary 
$\partial(W,S)$ has a cut pair (see remark before Theorem 2 in the introduction to [PS]), and hence it is not homeomorphic to $\Pi$.
This completes the proof.

\medskip
In view of the characterization of the Sierpi\'nski curve recalled above,
to prove the converse implication in Conjecture 0.1 (and in Theorem 2), it is sufficient
to show that $\partial(W,S)$ satisfies the properties from this characterization. Obviously, the boundary of any Coxeter system
is a compact metric space. Next three lemmas show that planarity
and unseparability of the labelled nerve $L^\bullet$ implies that $\partial(W,S)$
is connected, planar and 1-dimensional, respectively.

\medskip\noindent
{\bf 1.2 Lemma.}
{\it If $L^\bullet$ is planar,  unsepareable and distinct from a simplex then
$W$ is 1-ended and $\partial(W,S)$ is connected.}

\medskip\noindent
{\bf Proof:}
For a CAT(0) group connectedness of its visual boundary is
a consequence of its 1-endedness. Thus, in view of Theorem 8.7.2 in [Da],
connectedness of $\partial(W,S)$ follows from the fact that
$L$ is connected, has no separating simplex, and is distinct from a simplex.

\medskip\noindent
{\bf 1.3 Lemma.}
{\it If the nerve $L$ is planar then the boundary
$\partial(W,S)$ is also planar.}

\medskip\noindent
{\bf Proof:}
Extend $L$ to a triangulation $N$ of $S^2$ so that the following two
conditions hold:
\item{(1)} $L$ is a full subcomplex of $N$;
\item{(2)} $N$ is {\it flag relative to} $L$, i.e. if $T$ is a set of vertices of $N$
which are pairwise connected with edges, and if
$T\cap L$ spans a simplex of $L$, then $T$
spans a simplex of $N$.

\noindent
(We omit an elementary argument showing that such an extension is always possible for a planar complex.) Consider the labelled complex $N^\bullet$ whose labelling coincides with the labelling of $L^\bullet$
on edges contained in $L$, and which has labels equal to 2 at all
other edges of $N$. 
It is not hard to observe that $N^\bullet$ is then the labelled nerve of
some Coxeter system, which we denote 
$(W_{N^\bullet},S_{N^\bullet})$.
Moreover, $W$ is then a special subgroup of $W_{N^\bullet}$,
with the induced labelled nerve equal to $L^\bullet$.
Consequently, by Theorem A.1 from the appendix,
$\partial(W,S)\subset\partial(W_{N^\bullet},S_{N^\bullet})$.
Since the latter is homeomorphic to the sphere $S^2$
(see e.g. Theorem 3b.2 in [DJ] or Corollary 1 in [Dr]),
it follows that $\partial(W,S)$ is planar.

\medskip\noindent
{\bf 1.4 Lemma.}
{\it If $L^\bullet$ is planar,  unseparable, distinct from a simplex
and from a triangulation of the 2-sphere $S^2$, then
$\dim\partial(W,S)=1$.}

\medskip\noindent
{\bf Proof:}
Denote by $\hbox{vcd}(W)$ the virtual cohomological dimension of $W$.
It follows from results of Mike Davis that
$$
\hbox{vcd}(W)=\max\{ n\,:\,\overline{H}^{n-1}(L\setminus\sigma)
\ne0, \hbox{ for some simplex $\sigma$ of $L$, 
%including the empty simplex
or }\overline{H}^{n-1}(L)
\ne0  \}
$$
(see Corollary 8.5.5 in [Da]). Moreover, by the fact that $W$ acts geometrically
on $\Sigma$ and is virtually torsion-free, the pair 
$(\Sigma\cup\partial(W,S),\partial(W,S))$ is a ${\cal Z}$-structure
(in the sense of Bestvina described in [Be]) for any
torsion-free finite index subgroup $H<W$.
Since, by Theorem 1.7 of [Be] we then have $\dim\partial(W,S)=\hbox{cd}(H)-1$
(where $\hbox{cd}(H)$ is the cohomological dimension of $H$),
it follows that $\dim\partial(W,S)=\hbox{vcd}(W)-1$. Since $L$
is planar and distinct from $S^2$, we get 
(from the above formula of Davis) that
$\hbox{vcd}(W)\le2$,
and hence $\dim\partial(W,S)\le1$. On the other hand,
by Corollary 8.5.6 in [Da], a Coxeter group $W$
is virtually free iff $\dim\partial(W,S)\le0$. Since under our assumptions
$W$ is 1-ended (see Lemma 1.2), we conclude that 
$\dim\partial(W,S)=1$, thus completing the proof.

\bigskip\noindent
{\bf 2. Hyperbolic Coxeter groups and completion of the proof of Theorem 2.}

\medskip
In this section we assume that the group $W$ in a Coxeter system $(W,S)$
is word hyperbolic. Under this assumption, we deal with the conditions
of local connectedness and of absence of local cutpoints in the boundary
$\partial W=\partial(W,S)$.
 This allows to conclude the proof of Theorem 2.

\medskip\noindent
{\bf 2.1 Lemma.}
{\it If $L^\bullet$ is planar,  unseparable, distinct from a simplex,
from a triangulation of $S^1$ and $S^2$,
and from a labelled wheel,
and if $W$ is word-hyperbolic, then
the boundary $\partial W=\partial(W,S)$ is locally connected
and has no local cut point.}

\medskip\noindent
{\bf Proof:} 
%If $W$ is word-hyperbolic, its boundary $\partial(W,S)$ coincides
%with the Gromov boundary $\partial W$. 
By a result of Bestvina and Mess
[BeM], Gromov boundary of any 1-ended word-hyperbolic group
is locally connected. Hence, by Lemma 1.2, $\partial(W,S)$
is locally connected.

Unseparability of $L^\bullet$ means exactly 
that $W$ does not {\it visually} split 
(in the sense of the paper [MT] by Mihalik and Tschantz)
over a finite or 2-ended subgroup.
More precisely, this means that $W$ cannot be expressed
as an essential free product of its two special subgroups,
amalgamated along a finite or 2-ended special subgroup.
It follows from the main result of the same paper [MT]
that, under this condition, $W$ does not split along any finite or 2-ended subgroup.
Then, by a result of Bowditch [B], the boundary
$\partial W=\partial(W,S)$ has no local cut point,
or $W$ is a cocompact Fuchsian group
(i.e. a group of isometries of the hyperbolic plane acting
properly and cocompactly).
The latter case can be excluded as follows.
By a result of Davis
(see Theorem B in [D1] or Theorem 10.9.2 in [Da]),
if $W$ is a cocompact Fuchsian group then its nerve 
is either a triangulation of $S^1$ or 
$W$ splits as the direct sum of a special subgroup with the nerve $S^1$,
and another special subgroup, which is finite.
Since these possibilities are clearly inconsistent with the assumptions
of the lemma, the proof is completed.

\medskip\noindent
{\bf Proof of Theorem 2:}
In view of the Whyburn's characterization of the Sierpi\'nski curve
(recalled at the beginning of Section 1),
%comment preceding Lemma 1.2, 
Theorem 2 follows from
Lemmas 1.1-1.4 and 2.1.

\medskip\noindent
{\bf 2.2 Remarks.}
\item{(1)}
Arguments as in the proof of Lemma 2.1 cannot be easily extended to cover
the case of Coxeter groups that are not word-hyperbolic. For example,
it is not known fully for which 1-ended Coxeter groups $W$ the boundary
$\partial(W,S)$ is locally connected. Some criteria in the right angled case
are provided in [MRT] and [CM]. It is quite possible that using these criteria
one can show that for right angled Coxeter groups $W$, under assumptions
on $L$ as in Lemma 2.1, the boundary $\partial(W,S)$ is always locally connected.
Anyway, no criteria for local connectedness of $\partial(W,S)$ are known for
$W$ which are neither word-hyperbolic nor right angled.

\item{(2)} If $W$ is not word-hyperbolic, nothing seems to be known
about non-appearance of local cut points in the boundary 
$\partial(W,S)$. The main result of the paper [PS] by Papasoglu and Swenson allows,
under our assumptions on the labelled link $L^\bullet$,  
to exclude the appearance of cut pairs in $\partial(W,S)$. We don't know
if the methods developed in [PS] could be applied to exclude
as well local cut points (even assuming that $\partial(W,S)$ is locally
connected).

\bigskip

\noindent
{\bf Appendix.} %[January 20, 2015]
%boundaries of special subgroups embed in the boundary
%of a Coxeter group.

\medskip
In this appendix we provide a proof of the following folklore
result, which seems to have no explicitely presented proof in the literature.

\medskip\noindent
{\bf A.1 Proposition.}
{\it Let $(W,S)$ be any Coxeter system, and let $W_T$ be the special subgroup of $W$ corresponding to a subset $T\subset S$. Then the boundary
$\partial(W_T,T)$ embeds as a subspace in the boundary
$\partial(W,S)$.}

\medskip
Our proof of Proposition A.1 goes backwards. It consists of a sequence
of reductions to some other facts of independent interest, namely
Proposition A.2, Lemma A.3 and Lemma A.8. These reductions are
fairly standard and well known. The final step consists of proving
Lemma A.8, which we do by borrowing some arguments from
Ian Leary's paper [L] (where they were used for slightly different
purposes).

Recall that to each Coxeter system $(W,S)$ there is associated its
Coxeter-Davis complex $\Sigma(W,S)$. This complex is equipped with the piecewise euclidean Moussong's metric,  for which it is a CAT(0) space.  
The boundary
$\partial(W,S)$ is by definition the visual boundary of $\Sigma(W,S)$.
%equipped with this metric. 
Recall also that the Coxeter-Davis complex
$\Sigma(W_T,T)$ of a special subgroup is canonically a subcomplex
in $\Sigma(W,S)$ (see Section 7.3, especially Proposition 7.3.4, in [Da]).
Consequently,
since the visual boundary $\partial Y$ of a convex subspace $Y$ in a CAT(0) space $X$
canonically embeds in the visual boundary $\partial X$, to prove Proposition A.1,
it is sufficient to prove the following.

\medskip\noindent
{\bf A.2 Proposition.}
{\it Given a Coxeter system $(W,S)$, and a special subgroup 
$W_T<W$, the subcomplex $\Sigma(W_T,T)$ is a convex subspace
in the Coxeter-Davis complex $\Sigma(W,S)$ (equipped with the Moussong's metric).}

\medskip
To prove Proposition A.2 we need some preparations.
Recall that a subspace $K$ of a $CAT(1)$ space $L$ is $\pi$-convex
if any geodesic in $L$ of length less than $\pi$, having both endpoints
in $K$, is entirely contained in $K$. Recall also that if $X$ is a CAT(0)
piecewise euclidean complex then its any vertex link is a piecewise
spherical CAT(1) complex. Moreover, a connected subcomplex $Y$
in a CAT(0) piecewise euclidean complex $X$ is convex iff
for any vertex $v$ of $Y$ the link of $Y$ at $v$ is a $\pi$-convex
subspace in the link of $X$ at $v$. 
Since all vertex links in a Coxeter-Davis complex $\Sigma(W,S)$
are isometric to the nerve $L(W,S)$ of the system $(W,S)$,
to prove Proposition A.2, it is sufficient to prove the following.

\medskip\noindent
{\bf A.3 Lemma.}
{\it Given a Coxeter system $(W,S)$, and a special subgroup 
$W_T<W$, the nerve $L(W_T,T)$, viewed canonically as
the subcomplex in the nerve $L(W,S)$, is a $\pi$-convex subspace.}

\medskip
To prove Lemma A.3,
we will need the following terminology and tools introduced by 
Gabor Moussong.

\medskip\noindent
{\bf A.4 Definition.}
{A spherical simplex $\sigma$ is {\it of size $\ge{\pi\over2}$}
if its every edge has length $\ge{\pi\over2}$.}

\medskip\noindent
{\bf A.5 Defintion} (cf. Definition I.7.1 and Lemma I.5.1 in Appendix I in [Da]). 
A piecewise spherical simplicial complex $L$, with all cells of
size $\ge{\pi\over2}$, is {\it metrically flag}
if any of the following two equivalent conditions holds:
\item{$\bullet$} if $v_0,\dots,v_k$ are vertices of $L$ pairwise connected with edges, and $c_{ij}=d_L(v_i,v_j)$, then $v_0,\dots,v_k$ span
a simplex of $L$ iff the matrix $(c_{ij})$ is positive definite;
\item{$\bullet$} if $v_0,\dots,v_k$ are vertices of $L$ pairwise connected with edges, and $c_{ij}=d_L(v_i,v_j)$, then $v_0,\dots,v_k$ span
a simplex of $L$ iff $(c_{ij})$ is the matrix of edge lengths of
a spherical $k$-simplex.

\medskip\noindent
{\bf A.6 Lemma} (cf. Lemma 12.3.1 in [Da]). 
{\it For any Coxeter system $(W,S)$ its nerve $L(W,S)$, with its
natural piecewise spherical structure, has all simplices of size $\ge{\pi\over2}$,
and is metrically flag.}

%\medskip
%Before proving Lemma A.7, we recall one more fact due to Moussong.

\medskip\noindent
{\bf A.7 
%Moussong's 
Lemma} (cf. Lemma I.7.4 in Appendix I in [Da]).
{\it Let $L$ be a finite piecewise spherical simplicial complex with all cells
of size $\ge{\pi\over2}$. Then $L$ is $CAT(1)$ iff it is metrically flag.}

\medskip
In view of Lemma A.6, and since $L(W_T,T)$ is obviously
a full subcomplex of $L(W,S)$, to prove Lemma A.3,
it is sufficient to prove the following.

\medskip\noindent
{\bf A.8 Lemma.}
{\it Let $L$ be a finite piecewise spherical simplicial complex with all cells
of size $\ge{\pi\over2}$, and suppose it is CAT(1). Then its any full
subcomplex is a $\pi$-convex subspace of $L$.}

\medskip\noindent
{\bf Proof:}
We adapt the arguments from the proof of Theorem B.7 in Appendix B
of the paper [L] by Ian Leary.

Note that, since $L$ is CAT(1), it follows from Lemma A.7 that it is
metrically flag. 
Consider the simplicial complex $L*_KL$ obtained by gluing the two copies
of $L$ identically along the subcomplex $K$. Equip it with the piecewise
spherical metric induced from the metrics of the two copies of $L$.
Observe that $L*_KL$ has then cells of size $\ge{\pi\over2}$.
Moreover, since $K$ is a full subcomplex,
it is easy to see that $L*_KL$ is metrically flag.
Consequently, 
applying again Lemma A.7, it is CAT(1).

We claim that the metric of $L*_KL$ restricted to any copy of $L$ coincides with the original metric of $L$. One inequality between these metrics is obvious, because $L$ is a subcomplex of $L*_KL$. The other inequality follows easily from existence of the isometric involutive automorphism exchanging the two copies of $L$ in $L*_KL$.

Let $\gamma$ be a geodesic of length less than $\pi$ in $L$, with both
endpoints in $K$. By the above claim, lift of $\gamma$ to the first copy of $L$ in $L*_KL$ is also a geodesic; denote it by $\gamma_1$. If $\gamma$ is not entirely contained in $K$, the image of $\gamma_1$ through the above mentioned involution exchganging the two copies of $L$ is a different geodesic in $L*_KL$ connecting the same endpoints. This contradicts uniqueness of geodesics of length less than $\pi$ in $CAT(1)$ spaces,
thus completing the proof.

\bigskip

\bigskip\noindent
{\bf References}

\medskip

\item{[Be]} M. Bestvina, {\it Local homology properties of boundaries of groups,}
Michigan Math, J. 43 (1996), 123--139.

\item{[BeM]} M. Bestvina, G. Mess, {\it The boundary of negatively curved
groups}, J. Amer. Math. Soc. 4 (1991), 469--481.

\item{[B]} B. Bowditch, {\it Cut points and canonical splittings
of hyperbolic groups},  Acta Math. 180 (1998), 145--186.

\item{[BH]} M. Bridson, A. Haefliger, 
Metric Spaces of Non-Positive
Curvature, Grundlehren der mathematischen Wissenschaften 319, Springer,
1999.

\item{[CM]} W. Camp, M. Mihalik,
{\it A classification of right-angled Coxeter groups with no 3-flats and locally connected boundary}, J. Group Th. 17 (2014), 717--756.

\item{[Da]} M. Davis, {The geometry and topology 
of Coxeter groups}, {London Mathematical Society Monographs Series},
{vol. 32}, {Princeton University Press}, {Princeton}, {2008}.

\item{[D1]} M. Davis, {\it The cohomology of a Coxeter group with group
ring coefficients}, Duke Math. J. 91 (1998), 297--314.

\item{[DJ]} M. Davis, T. Januszkiewicz, {\it Hyperbolization of polyhedra},
J. Diff. Geom. 34 (1991), 347--388.

\item{[Dr]} A. Dranishnikov, {\it On boundaries of hyperbolic Coxeter groups}, Topology and its Applications 110 (2001), 29--38. 

\item{[KB]}
I. Kapovich, N. Benakli, {\it Boundaries of hyperbolic groups,} 
in: Combinatorial and geometric group theory 
(New York, 2000/Hoboken, NJ, 2001), 39--93,
Contemp. Math. 296, Amer. Math. Soc., Providence, RI, 2002.

\item{[L]} I. Leary, {\it A metric Kan-Thurston theorem},
Journal of Topology 6 (2013), 251--284. %doi:10.1112/jtopol/jts035.

\item{[MT]} M. Mihalik, S. Tschantz, {\it Visual decompositions
of Coxeter groups}, Groups Geometry and Dynamics 3 (2009), 173--198. 

\item{[MRT]} M. Mihalik, K. Ruane, S. Tschantz, {\it Local connectivity
of right-angled Coxeter group boundaries}, J. Group Theory 10 (2007),
531--560.

\item{[PS]} P. Papasoglu, E. Swenson, {\it Boundaries and 
JSJ decompositions of CAT(0) groups}, Geom. Funct. Anal. 19 (2009),
558--590. %DOI 10.1007/s00039-009-0012-8. 

\item{[Wh]} G. T. Whyburn, {\it Topological characterization of the 
Sierpinski curve}, Fundamenta Math. 45 (1958), 320--324.

\bye